\documentclass[reqno,12pt]{amsart}
\usepackage{amsfonts,amssymb,verbatim}
\usepackage{enumerate}

\newcommand{\Z}{\mathbb{Z}}

\newtheorem{thm}{Theorem}
\newtheorem{lem}{Lemma}

\newtheorem{cor}{Corollary}
\theoremstyle{definition}

\theoremstyle{remark}

\numberwithin{equation}{section}

\author{Jing-Jing Huang}
\address{
Jing-Jing Huang: Department of Mathematics, University of Toronto, Ontario M5S 2E4, Canada}
\email{huang@math.toronto.edu}
\thanks{Research is partially supported by the NSERC Grant A5123.}

\begin{document}

\title
[Rational points near planar curves]
{Rational points near planar curves and Diophantine approximation}

\begin{abstract}
In this paper, we establish asymptotic formulae with optimal errors for the number of rational points that are close to a planar curve, which unify and extend the results of Beresnevich-Dickinson-Velani \cite{BDV} and Vaughan-Velani \cite{VV}.
Furthermore, we complete the Lebesgue theory of Diophantine approximation on weakly non-degenerate planar curves that was initially developed by Beresnevich-Zorin \cite{BZ} in the divergence case. 
\end{abstract}
\maketitle

\section{Introduction} \label{s1}

The problem of counting rational points lying near manifolds is of fundamental interest in number theory. Both the arithmetic and geometric properties of the manifold come into play. More importantly, this problem is related to central theme in the metric theory of Diophantine approximation on manifolds, namely the \emph{Generalized Khintchine-Jarnik type problems}. In this paper we mainly study the case of planar curves $\mathcal{C}$ embedded in $\mathbb{R}^2$. This model case has been studied for a long time but was not reasonably understood until recently. We will discuss more about the nature of the associated counting problem and its applications in \S{\ref{s5}}, \S{\ref{s6}} and \S{\ref{s7}}. 

\subsection{Rational points near curves }
To motivate the problem, we first introduce some notations. First of all, as is well understood (see \cite{BD} for instance), we may assume our curve is locally represented as the graph $\mathcal{C}_f$ of some real valued function $f$ on a compact interval $I:=[\eta,\xi]\subset\mathbb{R}$, where $\mathcal{C}_f=\{(x,f(x))|x\in I\}$. With metrical applications in mind, this is certainly without loss of generality.

With the goal of counting rational points near the curve $\mathcal{C}_f$, we define, as in \cite{Hu,BDV,VV,BZ}, for any $1/2\ge\delta>0$ and $Q\ge1$,
$$\widetilde{N}_f(Q,\delta):=\#\left\{(a/q,b/q)\in\mathbb{Q}^2:\begin{array}{l}
q\le Q, a/q\in I, (a,b,q)=1\\
|f(a/q)-b/q|<\delta/Q
\end{array}
\right\}.$$
 Roughly speaking, this function counts the number of rational points with denominator at most $Q$ which lie within distance $\delta Q^{-1}$ of the curve $\mathcal{C}_f$.  

If $f$ is linear (namely $C_f$ is an affine line),  the quantity $\widetilde{N}_f(Q,\delta)$ depends crucially on how well the slope can be approximated and hence is well understood (see \cite{Sp} and \cite{BD} for more detailed discussions). Therefore, it is natural to consider the following class of functions with non-vanishing curvature.
Suppose $c_1$ and $c_2$ are two fixed constants with $0<c_1\le c_2$. Now we denote by $\mathcal{F}(I;c_1,c_2)$ the set of $C^2$ functions $f: I\rightarrow\mathbb{R}$ such that
$$c_1\le|f''(x)|\le c_2\;\;\;\text{ for all }x\in I.$$

Huxley \cite{Hu} established, via elementary geometry and Swinnerton-Dyer's ingenious determinant method \cite{Sw}, that for any $f\in\mathcal{F}(I;c_1,c_2)$, one has
\begin{equation}\label{e1.0}
\widetilde{N}_f(Q,\delta)\ll \delta^{1-\varepsilon}Q^2+Q
\end{equation}
where the implicit constant only depends on $c_1$, $c_2$, $\varepsilon$ and the length of the interval $I$. We note that Huxley works on the range $Q\le q<2Q$ but it is trivially seen that summing over the dyadic ranges gives essentially the same bound for the range $q\le Q$.

Using analytic methods, Vaughan and Velani \cite{VV} remarkably removed the $\varepsilon$ in Huxley's bound. They showed that for any $f\in\mathcal{F}(I;c_1,c_2)$
\begin{equation}\label{e1.2}
\widetilde{N}_f(Q,\delta)\ll \delta Q^2+\delta^{-\frac12}Q
\end{equation}
and that, with the additional assumption $f''\in \text{Lip}_\theta(I)$ with $\theta\in(0,1)$, for any $\varepsilon>0$
\begin{equation}\label{e1.3}
\widetilde{N}_f(Q,\delta)\ll \delta Q^2+\delta^{-\frac12}Q^{\frac12+\varepsilon}+\delta^{\frac{\theta-1}2}Q^\frac{3-\theta}2.
\end{equation}

Based on the seminal work of Kleinbock and Margulis \cite{KM} on extremal manifolds and using the ubiquitous system developed by themselves,  Beresnevich, Dickinson and Velani established in a breakthrough paper \cite{BDV} a sharp lower bound in the range $\delta\gg Q^{-1}$. More precisely, they showed that for any $f\in\mathcal{F}(I;c_1,c_2)$ satisfying the additional condition $f\in C^3(I)$ there exist constants $k_1,k_2,c,Q_0$ such that for any $Q>Q_0$ and any $\delta$ with $k_1Q^{-1}\le\delta\le k_2$ one has

\begin{equation}\label{e1.4}
\widetilde{N}_f(Q,\delta)\ge c\delta Q^2.
\end{equation}

Note however that the constant $c$ above is not obtained explicitly and that it depends on $f$.
Later on, Beresnevich and Zorin \cite{BZ} were able to get an explicit constant $c$ in \eqref{e1.4} in terms of $c_1$ and $c_2$ and also extend \eqref{e1.4} to a topologically closed set of functions.

The set $\mathcal{F}(I;c_1,c_2)$ is clearly a subset of the Banach space $\mathcal{C}(I)$ of continuous functions on $I$ with the uniform norm. 
Now let $\overline{\mathcal{F}}(I;c_1,c_2)$ be the topological closure of $\mathcal{F}(I;c_1,c_2)$ under the uniform convergence norm that is equipped on $\mathcal{C}(I)$. 

It is shown in \cite{BZ} that \eqref{e1.4} holds for any $f\in\overline{\mathcal{F}}(I;c_1,c_2)$ with the explicit constant
\begin{equation}\label{e1.20}
c=\frac{{c_1}^{13}(\min\{1,\sqrt{c_1}\})^{12}}{2^{89}3^{72}{c_2}^{15}}|I|.
\end{equation}
The above constant $c$  is rather small and certainly far from the expected value, while implicit constants in the upper bounds \eqref{e1.2} and \eqref{e1.3} seem to have never been worked out. One can in principle trace back in the argument of \cite{VV} and obtain some explicit constants, which are unfortunately doomed not to be optimal.

In this paper, we first of all work with a slight variant of the above counting function $\widetilde{N}_f(Q,\delta)$. Let
\begin{equation}\label{e1.1}
\hat{N}_f(Q,\delta):=\#\left\{\left(\frac{a}q,\frac{b}q\right)\in\mathbb{Q}^2:\begin{array}{l}
q\le Q, a/q\in I, (a,b,q)=1\\
\left|f\left(\frac{a}q\right)-\frac{b}q\right|<\frac\delta{q}
\end{array}
\right\}.
\end{equation}

In one of our main theorems, we establish, for the first time, an asymptotic formula for $\hat{N}_f(Q,\delta)$ in the range $\delta>Q^{-1+o(1)}$, which immediately implies good upper and lower bounds for $\widetilde{N}_f(Q,\delta)$ (see Theorem \ref{t7} below).

Let $\zeta(s)$ be the Riemann zeta-function.

\begin{thm}\label{t1}
Let $f\in{\mathcal{F}}(I;c_1,c_2)$ have a Lipschitz continuous second derivative. Then, for any $\varepsilon>0$, $1/2\ge\delta>0$ and integer $Q>1$, we have
$$\hat{N}_f(Q,\delta)=\frac1{\zeta(3)}|I|\delta Q^2+O\left(\delta^{\frac12}\left(\log \frac1\delta\right)Q^{\frac32}+Q^{1+\varepsilon}\right)$$
where the implicit constant only depends on $I$, $c_1$, $c_2$, $\varepsilon$ and the Lipschitz constant, and is in particular independent of $f$, $\delta$ and $Q$. 
\end{thm}

It should be pointed out that the bound Theorem \ref{t1} renders is almost best possible. To see that, we make the following observations.

 The main term in Theorem \ref{t1} represents the expected number of rational points satisfying \eqref{e1.1}, and it is in particular independent on $c_1$ and $c_2$ and uniform for all $C^3$ curves. Moreover, it is completely in compliance with the following probabilistic heuristics. Consider the set $\mathcal{G}_{I,Q}$ of  fractions $a/q\in I$ with $q\le Q$ ($a$ and $q$ are not necessarily coprime). For each such $a/q$, the possibility that there exists $b$ such that 
$$\left|f\left(\frac{a}q\right)-\frac{b}q\right|<\frac\delta{q}$$
is $2\delta$ if we assume $\|qf(a/q)\|$ behaves like a random number in $[0,1/2]$. Also there are in total roughly 
$$|I|\sum_{q\le Q}q\sim \frac1{2}|I|Q^2$$
 fractions $a/q$ in $\mathcal{G}_{I,Q}$. Moreover, it is an easy exercise to show that the probability that $(a,b,q)=1$ is $\zeta(3)^{-1}$ again assuming $a$, $b$ and $q$ are independent random integers. Hence this probabilistic model suggests that 
$$\hat{N}_f(Q,\delta)\sim \frac1{\zeta(3)}|I|\delta Q^2,$$
which coincides with the main term in Theorem \ref{t1}.
However, we will soon see below that this is not true unconditionally.

Our main term dominates only when $\delta\gg Q^{-1+o(1)}$. In the opposite case, $\delta\ll Q^{-1+o(1)}$, the error $Q^{1+\varepsilon}$ dominates. Moreover $\delta^{\frac12}\left(\log \frac1\delta\right)Q^{\frac32}$ is bounded by the geometric mean of $\delta Q^2$ and $Q^{1+\varepsilon}$ when $\delta\ge Q^{-1}$ and by $Q^{1+\varepsilon}$ when $\delta<Q^{-1}$. So
$$\delta^{\frac12}\left(\log \frac1\delta\right)Q^{\frac32}\ll\max\{\delta Q^2, Q^{1+\varepsilon}\}.$$
Now we consider the example $f(x)=x^2, x\in[1,2]$ to show that the error $O(Q^{1+\varepsilon})$ is tight. All the points $\{a/q, a^2/q^2\}$ with $q\le \sqrt{Q}$ and $q\le a\le 2q$, $(a,q)=1$ lie on the parabola $f(x)=x^2$ and are clearly counted in $\hat{N}_f(Q,\delta)$ for any $\delta>0$. There are roughly $\phi(q)$ of them for each denominator $q$. This yields
\begin{equation}\label{e1.5}
\hat{N}_f(Q,\delta)\gg \sum_{q\le\sqrt{Q}}\phi(q)\gg Q
\end{equation}
regardless of the value of $\delta$. But Theorem \ref{t1} implies that
$$\hat{N}_f(Q,\delta)\ll Q^{1+\varepsilon}$$
as $\delta$ goes to 0, 
which explains why the exponent in $Q^{1+\varepsilon}$ is not subject to improvement except for the $\varepsilon$ power.

Theorem \ref{t1} immediately implies the following corollary.

\begin{cor}\label{c1}
Under the same conditions as in Theorem \ref{t1},  when $1/2\ge\delta>Q^{-1+\varepsilon}$ for some $\varepsilon\in(0,1)$ we have
$$\hat{N}_f(Q,\delta)=\left(\frac1{\zeta(3)}+O(Q^{-\frac\varepsilon2})\right)|I|\delta Q^2.$$
\end{cor}

A consequence of this is the following result on $\widetilde{N}_f(Q,\delta)$, the proof of which is presented in \S{\ref{s1.5}}.

\begin{thm}\label{t7}
Under the same conditions as in Theorem \ref{t1},  when $1/2\ge\delta>Q^{-1+\varepsilon}$ for some $\varepsilon\in(0,1)$  we have
$$\frac{2\sqrt{3}}{9\zeta(3)}+O(Q^{-\frac\varepsilon2})\le\frac{\widetilde{N}_f(Q,\delta)}{|I|\delta Q^2}\le \frac1{\zeta(3)}+O(Q^{-\frac\varepsilon2}).$$
\end{thm}

It is possible that an adaptation of our method would give the asymptotics for $\widetilde{N}_f(Q,\delta)$. We expect that 
$${\widetilde{N}_f(Q,\delta)}\sim\frac2{3\zeta(3)}{|I|\delta Q^2}$$
holds for the same range of $\delta$ and $Q$ as above.
In any event, the constant $\frac{2\sqrt{3}}{9\zeta(3)}$ in our lower bound should be compared with that in \eqref{e1.20}. In particular, neither of the leading constants in our upper and lower bounds for $\widetilde{N}_f(Q,\delta)$ depends on $c_1, c_2$, as opposed to the lower bound constant \eqref{e1.20} obtained by Beresnevich and Zorin \cite{BZ}.

Also, if $\delta<Q^{-1}$, one can no longer expect a nontrivial lower bound for the counting function ${\widetilde{N}_f(Q,\delta)}$ (as well as for other similar counting functions $\hat{N}_f(Q,\delta)$ and ${N}_f(Q,\delta)$ to be defined below), which depends very much on the arithmetic property of the specific curve under consideration. We will discuss this phenomenon in \S{\ref{s5}}.

To prove Theorem \ref{t1}, it is more convenient to work without the coprime condition $(a,b,q)=1$ then we may recover the condition by using the M\"{o}bius inversion. Let
$$
N_f(Q,\delta):=\#\left\{(a,q)\in\mathbb{Z}\times\mathbb{N}:\begin{array}{l}
q\le Q, a/q\in I\\
\|qf(a/q)\|<\delta
\end{array}
\right\}
$$
where we use the standard notation $\|x\|=\inf_{k\in \Z}|x-k|$. Note that $a$ and $q$ are not required to be coprime in the above definition. This form is particularly suited for the application to Diophantine approximation.

\begin{thm}\footnote{After the submission of the present work, it was brought to the author's attention that Theorem \ref{t6} with a somewhat different error term had been recently proved independently by A. Gafni \cite{Ga}. }\label{t6}
Under the same conditions as in Theorem \ref{t1}, we have
$$N_f(Q,\delta)=|I|\delta Q^2+O\left(\delta^{\frac12}\left(\log \frac1\delta\right)Q^{\frac32}+Q^{1+\varepsilon}\right)$$
where the implicit constant only depends on $I$, $c_1$, $c_2$, $\varepsilon$ and the Lipschitz constant, and is in particular independent of $f$, $\delta$ and $Q$.
\end{thm}

Without the Lipschitz condition in Theorem \ref{t1}, we are less successful. Nevertheless, we can prove the following weaker asymptotics for functions in $\overline{\mathcal{F}}(I;c_1,c_2)$.

\begin{thm}\label{t2}
Let $f\in\overline{\mathcal{F}}(I;c_1,c_2)$. Then for any $\delta\in(0,\frac12]$ and integer $Q>1$ we have
$$N_f(Q,\delta)=|I|\delta Q^2+O\left(\delta^{\frac12}\left(\log \frac1\delta\right)Q^{\frac32}+Q^{\frac43}\right)$$
where the implicit constant only depends on $I$, $c_1$ and $c_2$ and is in particular independent of $f$, $\delta$ and $Q$.
\end{thm}

One can make a similar discussion about the sizes of the main term and the error. This time the main term wins when $\delta\gg Q^{-2/3}$. However, it remains interesting to see whether our error is optimal in this more general context.

We remark in passing that the set $\overline{\mathcal{F}}(I;c_1,c_2)$ is genuinely larger than $\mathcal{F}(I;c_1,c_2)$ in that it is shown in \cite{BZ} that there are functions in $\overline{\mathcal{F}}(I;c_1,c_2)$ whose second derivatives do not exist in a dense subset of $I$. It is also shown in \cite{BZ} that $\overline{\mathcal{F}}(I;c_1,c_2)$ is the same as the set of differentiable functions on $I$ whose derivatives are bi-Lipschitz with constants $c_1$ and $c_2$.



\subsection{Metric Diophantine approximations}

There has been an explosive growth of the literature on the subject of metric Diophantine approximation during the past few decades, hence we do not intend to give a comprehensive reference here. Instead, we mention only a few results that are of direct relevance to the topics studied in this paper and refer the readers to the excellent book \cite{Sp} of Sprind\v{z}uk for a general introduction to this subject and the book of Bernik and Dodson \cite{BD} or the research papers \cite {BDV}, \cite{Be} for more recent developments.

The foundations of the modern theory of metrical Diophantine approximation on planar curves were laid by Schmidt in \cite{Sc1} and \cite{Sc2} in 1964, in which he proved that every $C^3$ non-degenerate planar curve is extremal. In particular, Schmidt's proof relies on a counting result which is reminiscent of but much weaker than our Theorems \ref{t6}, \ref{t2}. We recall that a curve $\mathcal{C}_f$ for some $C^2$ function $f: I\rightarrow\mathbb{R}$ is 
\emph{non-degenerate} if for almost all points on the curve the curvature does not vanish i.e. $f''(x)\not=0 \text{ a.e. }x\in I$. 

Let $\psi:\mathbb{N}\rightarrow\mathbb{R}^+$ be a test function. And let

$$\mathcal{S}_f(\psi):=\{ x\in I|\exists^\infty q\in\mathbb{N} \text{ such that } \max\{\|qx\|,\|qf(x)\|\}<\psi(q)\}$$
and
$$\mathcal{S}^*_f(\psi):=\{ x\in I|\exists^\infty q\in\mathbb{N} \text{ such that } \|qx\|\cdot\|qf(x)\|<\psi(q)\}$$

where $\exists^\infty$ stands for ``there exist infinity many$\cdots$".

The key feature of problems of this type is that there is always a \emph{Zero versus Full dichotomy} about the measure of the set in question, according as whether a certain volume series converges. In our case, these sets are exactly $\mathcal{S}_f(\psi)$ and $\mathcal{S}^*_f(\psi)$.
However, determining what class of functions(curves) satisfies this dichotomy was a long standing problem and over 40 years after Schmidt's result nothing was known apart from Bernik's Khintchine type theorem for convergence for the parabola \cite{Be1}. In the past decade, there has been however some tremendous progress on this problem and the general theory of simultaneous approximation on non-degenerate curves was developed by Beresnevich-Dickinson-Velani \cite{BDV} for the divergence theory and Vaughan-Velani \cite{VV} for the convergence theory. The distribution of rational points near curves plays a crucial role in their approaches. In particular, establishing sharp upper and lower bounds for $\widetilde{N}_f(Q,\delta)$ such as $\eqref{e1.2}$, $\eqref{e1.3}$ and $\eqref{e1.4}$ is the main obstacle of the problem. The theory was subsequently generalized to multiplicative/inhomogeneous settings \cite{BL,BVV,BV1,BV2}. With additional treatment, Beresnevich and Zorin \cite{BZ} extended the divergence part to the set of weakly non-degenerate planar curves.  

The curve $\mathcal{C}_f=\{(x,f(x))|x\in I\}$ will be called \emph{weakly non-degenerate at} $x_0\in I$ if there exist constants $c_2\ge c_1>0$ and a compact interval $I_0\subset I$ centered at $x_0$ such that $f|_{I_0}\in\overline{\mathcal{F}}(I_0;c_1,c_2)$. We say $\mathcal{C}_f$ is \emph{weakly non-degenerate} if $\mathcal{C}_f$ is weakly non-generate at almost every point $x_0\in I$. Geometrically, weakly non-degenerate curves are defined locally by functions in $\overline{\mathcal{F}}(I_0;c_1,c_2)$ except for a Lebesgue measure zero set.

\begin{thm}[Simultaneous approximation]\label{t3}
Let $\psi$ be a monotonic arithmetic function and $f:I\rightarrow\mathbb{R}$ be a weakly non-degenerate function. Then
\[
|S_f(\psi)|=\left\{
\begin{array}{lll}
\textnormal{\small{Z\tiny{ERO}}}&\textnormal{if}& \displaystyle\sum_{q=1}^{\infty}\psi(q)^2<\infty\\
\textnormal{\small{F\tiny{ULL}}}&\textnormal{if}&
\displaystyle\sum_{q=1}^{\infty}\psi(q)^2=\infty.
\end{array}
\right.
\]
\end{thm}
Here `full' means that $|S_f(\psi)|=|I|$. As mentioned above, the divergence part in this theorem was first established in \cite[Theorem 4(A)]{BZ}. The convergence part is new in this setting. Hence Theorem \ref{t3} completes the Lebesgue measure theory for weakly non-degenerate planar curves. The proof of the convergence part of Theorem \ref{t3} modulo Theorem \ref{t2} is routine, see \cite[\S6]{VV} for the outline of the proof. Similarly we also have

\begin{thm}[Multiplicative approximation]\label{t4}
Under the same assumptions as in Theorem \ref{t3}, we have
$$|S^*_f(\psi)|=
\textnormal{\small{Z\tiny{ERO}}}\quad\textnormal{if}\quad \displaystyle\sum_{q=1}^{\infty}\psi(q)\log q<\infty.$$
\end{thm}

Again the proof of this theorem modulo Theorem \ref{t2} is the same as that of Theorem 1 of \cite{BL}. One just needs to replace \cite[Theorem 1]{VV} used in \cite{BL}  by our Theorem \ref{t2}.
To the best of our knowledge, the divergence counterpart of Theorem \ref{t4} is still open, even in the case of non-degenerate curves. Using our Theorem \ref{t2} in combination with the ideas of \cite{BVV} it is straightforward to state and prove an inhomogeneous version of the above two theorems. Along the same lines of \cite[\S 7]{VV} one can use Theorem \ref{t2} to prove a Hausdorff measure version of Theorem \ref{t3} for weakly non-degenerate curves. But due to the shorter range of $\delta$ in Theorem \ref{t2}, the Jarnik type result obtained is not as sharp as that which is deduced from Theorem \ref{t1} for $C^3$ non-degenerate curves. We leave this to the reader.

All the implicit constants in $\ll$ and $O$ notations in this paper are allowed to depend on various fixed constants $|I|$, $c_1$, $c_2$, $C$ and $\varepsilon$. Nevertheless, they are independent of any varying parameter like $\delta$, $Q$ and $f$ considered in this paper.    

This paper is structured as follows. In \S{\ref{s1.5}} we show how to deduce Theorem \ref{t1} and \ref{t7} modulo Theorem \ref{t6}. In \S{\ref{s2}}, we prove Theorem \ref{t6} but postpone the proof of one of the key steps namely the method of stationary phase to \S{\ref{s4}}. In \S\ref{s3} we prove Theorem \ref{t2} by adapting the proof in \S\ref{s2}. We discuss a question raised by Barry Mazur that is closely related to the counting problem studied in this paper and give a general answer in \S{\ref{s5}}. One can also look at our main results from the perspective of uniform distribution; a variant of Theorem \ref{t1} in this flavor is presented in \S{\ref{s6}}. Finally, we discuss the general framework of our approach in \S{\ref{s7}} in particular the main difficulty that prevents us from obtaining sharp bounds for non-degenerate submanifolds of higher dimensional Euclidean spaces.

\section{Proof of Theorem \ref{t1} and \ref{t7}}\label{s1.5}

Here unlike $N_f(Q,\delta)$, every pair of rational points $(a/q,b/q)$ is counted exactly once in $\hat{N}_f(Q,\delta)$. By factoring out the greatest common divisor $d$ of $a, b, q$, one readily observes the relation
$$N_f(Q,\delta)=\sum_{d\le Q}\hat{N}_f(Q/d,\delta/d)$$
and that by M\"{o}bius inversion
$$\hat{N}_f(Q,\delta)=\sum_{d\le Q}\mu(d){N}_f(Q/d,\delta/d)$$
where $\mu(\cdot)$ is the M\"{o}bius function. By Theorem \ref{t6} we have
$$\hat{N}_f(Q,\delta)=\sum_{d=1}^\infty\frac{\mu(d)}{d^3}|I|\delta Q^2+O(\delta^{\frac12}\log(\delta^{-1})Q^{\frac32}+Q^{1+\varepsilon}).$$ On noticing $\sum_{d=1}^\infty\frac{\mu(d)}{d^3}=\zeta(3)^{-1}$,
we deduce Theorem \ref{t1} modulo Theorem \ref{t6}.

Regarding Theorem \ref{t7}, we make the following two observations:
$$\widetilde{N}_f(Q,\delta)\le \hat{N}_f(Q,\delta)$$ 
and for some fixed $\alpha\in(0,1)$
\begin{align*}
\widetilde{N}_f(Q,\delta)&\ge\#\left\{\left(\frac{a}q,\frac{b}q\right)\in\mathbb{Q}^2:\begin{array}{l}
\alpha Q\le q\le Q, a/q\in I, (a,b,q)=1\\
\left|f\left(\frac{a}q\right)-\frac{b}q\right|<\frac{\alpha\delta}{q}
\end{array}
\right\}\\
&=\hat{N}_f(Q,\alpha\delta)-\hat{N}_f(\alpha Q,\alpha\delta)\\
&\overset{Cor. \ref{c1}}{=}\left((\alpha-\alpha^3)\zeta(3)^{-1}+O(Q^{-\frac{\varepsilon}2})\right)|I|\delta Q^2.
\end{align*}
Now we optimize the leading coefficient by taking $\alpha=1/\sqrt{3}$ and Theorem \ref{t7} follows.

\section{Proof of Theorem \ref{t6}} \label{s2}
Our proof of Theorem \ref{t6} follows that of \cite[Theorem 3]{VV}; the new novel feature is the use of a more sophisticated Fourier weight.
The starting point of our approach is based on approximation of the characteristic function $\chi_{\scriptscriptstyle\Delta}(x)$ of the set $\Delta:=\{x\in\mathbb{R}:\|x\|<\delta\}$ with $\delta\le1/2$ by Fourier series. Though it would be natural to consider the exact Fourier series expansion of $\chi_{\scriptscriptstyle\Delta}(x)$, in practice the Gibbs phenomenon due to the jump discontinuity of $\chi_{\scriptscriptstyle\Delta}(x)$ brings additional difficulties. To avoid this issue, there are naturally two possible options: use a finite Fourier sum or use the infinite Fourier expansion of a smoothed version of $\chi_{\scriptscriptstyle\Delta}(x)$. In \cite{VV}, the authors take the former approach and indeed they use the Fej\'{e}r kernel. Here, we conveniently base our proof on a refinement of the Erd\H{o}s-Tur\'{a}n inequality in \cite[Chap. 1]{Mo}, the advantage of which is that it absorbs the Fourier coefficients and leaves a clean exponential sum to tackle. Other authors in the literature also use approximation functions constructed by I.M. Vinogradov, see \cite{BD}, \cite{Sp}, \cite{Vi}.

Let $u_n$ be a sequence of $N$ real numbers with $1\le n\le N$. Let $Z(N;\alpha,\beta)$ count the number of $n$ for which $u_n\in(\alpha,\beta)\pmod{1}$ with $\alpha<\beta<\alpha+1$. It is of great interest to study the discrepancy function 
\begin{equation}\label{e2.1}
D(N;\alpha,\beta):=Z(N;\alpha,\beta)-(\beta-\alpha)N.
\end{equation}

\begin{lem}\label{l1} 
For any positive integer $K$,
$$|D(N;\alpha,\beta)|\le\frac{N}{K+1}+2\sum_{k=1}^{K}b_k\left|\sum_{n=1}^Ne(ku_n)\right|$$
where
$$b_k=\frac1{K+1}+\min\left(\beta-\alpha,\frac1{\pi k}\right).$$
\end{lem}

A proof of this lemma can be found for instance in Chapter 1 of Montgomery \cite{Mo}. 

Now we apply Lemma \ref{l1} to the sequence \{$qf(a/q)$\},  $a/q\in I$ and $q\le Q$ with $\alpha=-\delta$ and $\beta=\delta$. Notice that here we do not require that $a/q$ is in reduced form and hence we are really counting those numbers with multiplicities as long as the denominators are bounded above by $Q$. The overall number of such rationals is 
$$\sum_{q\le Q}\sum_{a/q\in I}1=\sum_{q\le Q}(|I|q+O(1))=|I|Q^2/2+O(Q).$$
Therefore for any positive integer $K$,
\begin{align}
&N_f(Q,\delta)-|I|\delta Q^2\nonumber\\
\ll&\frac{Q^2}{K}+\delta Q+\sum_{k=1}^{K}b_k\left|\sum_{q\le Q}\sum_{a/q\in I}e(kqf(a/q))\right|\label{e2.0}
\end{align}
where
$$b_k=\frac1{K+1}+\min\left(2\delta,\frac1{\pi k}\right).$$
It will be seen in the proof that the optimal choice for $K$ is
\[
K=\max\left\{\left\lfloor\frac1\delta\right\rfloor,Q\right\}.
\]
For the time being, we assume 
$$\delta\ge Q^{-1}.$$ In this case, $K=Q$.

In order to understand the exponential sum on the right side of \eqref{e2.0}, we apply a version of the Poisson summation in \cite[Chap. 3, Theorem 8]{Mo}. 
\begin{lem}[The Truncated Poisson summation formula]\label{l2}
Let $g$ be a real-valued function, and suppose that $g'$ is continuous and increasing on $[c,d]$. Let $s=g'(c)$ and $t=g'(d)$. Then
$$\sum_{c\le n\le d}e(g(n))=\sum_{s-1\le j \le t+1}\int_c^de(g(x)-jx)dx+O(\log(2+t-s)).$$
\end{lem}
We put $J:=[\min f',\max f']$ and $J_k:=[k\min f'-1,k\max f'+1]$. Thus Lemma \ref{l2} implies when $f''>0$
\begin{equation}\label{e2.2}
\sum_{a\in qI}e(kqf(a/q))=\sum_{j\in J_k}\int_{qI}e(kqf(r/q)-jr)dr+O(\log(2+k|J|)).
\end{equation} 
Notice that $|J|\le c_2|I|$ and hence is bounded. Therefore the contribution of the error in \eqref{e2.2} to the right side of \eqref{e2.0} is
\begin{equation}\label{e2.14}
\ll Q\sum_{k=1}^{K}\left(\frac1{K+1}+\min\left(2\delta,\frac1{\pi k}\right)\right)\log k\ll Q(\log K)^2 
\end{equation}
By changing variable $r=qx$, the integral in \eqref{e2.2} becomes
\begin{equation}\label{e2.3}
q\int_{I}e(q(kf(x)-jx))dx.
\end{equation}
When $f''<0$, we apply Lemma \ref{l2} with $g(x)=-kqf(x/q)$ and then take the complex conjugate. It is easily seen that \eqref{e2.2} still holds in this case.

Now to evaluate the oscillatory integral \eqref{e2.3} precisely, we need to apply the method of stationary phase in harmonic analysis. Roughly speaking the idea is that the phase function oscillates slowly near a critical point and hence the integral will have a peak; when staying away from the critical points, the phase functions oscillates frequently which results in cancellation in the integral.

To avoid some unnecessary technicalities when there is a critical point very close to an end point of $I$, we shrink the interval $J_k$ a little bit to rule this case out. Let $\tilde{J}_k:=[k\min f'+1,k\max f'-1]$. In the rare case that $k\max f'-1<k\min f'+1$ we just take $\tilde{J}_k$ to be empty. Then by Lemma 2 in \cite[Chapter 3]{Mo}
\begin{equation}\label{e2.20}
\int_{I}e(q(kf(x)-jx))dx\ll \frac1{\sqrt{qk}}
\end{equation}
since 
$$c_1\le|f''(x)|\le c_2\;\;\;\text{ for all }x\in I.$$
Thus
\begin{align}
&\sum_{k\le K}b_k\left|\sum_{q\le Q}\sum_{j\in J_k\backslash\tilde{J}_k}q\int_{I}e(q(kf(x)-jx))dx\right|\nonumber\\
\ll&\sum_{k\le K}\left(\frac1{K+1}+\min\left(\delta,\frac1{ k}\right)\right)\frac1{\sqrt{k}}\sum_{q\le Q}\sqrt{q}\nonumber\\
\ll&K^{-1/2}Q^{3/2}+\delta^{1/2}Q^{3/2}.
\label{e2.4}
\end{align}

\begin{lem}[The method of stationary phase]\label{l3}Let $I=[a,b]$ be a compact interval.
Suppose $\phi(x)\in\mathcal{F}(I;c_1,c_2)$ with $0<c_1\le c_2$ has a Lipschitz second derivative, $|\phi''(x_1)-\phi''(x_2)|
\le C|x_1-x_2|$ for some constant $C>0$ and for any $x_1,x_2\in I$.  Suppose $\phi'(x_0)=0$ for some $x_0\in(a,b)$. Then when $\lambda\rightarrow\infty$,
$$\int_{I}e(\lambda\phi(x))dx=\frac1{\sqrt{\lambda|\phi''(x_0)|}}e\left(\lambda\phi(x_0)\pm\frac18\right)+O\left(\left(\frac1\kappa+\log\lambda\right)\lambda^{-1}\right)$$
where the sign of $1/8$ is chosen as that of $\phi''(x)$ and $\kappa=\min\{x_0-a,b-x_0\}$. 
\end{lem} 
We note that the implied constant in the above lemma only depends on $I$ and $c_1,c_2,C$ and is in particular independent of $\phi$. The proof of this lemma is postponed to \S{\ref{s4}}. The novel feature here is that we put minimum differentiability condition on the phase function $\phi$, as opposed to many standard reference books (see \S{\ref{s4}} for more discussions).

Now we apply Lemma \ref{l3} with, $I=[\eta,\xi]$,   $\phi(x)=f(x)-j/k$ and $\lambda=qk$. For each $j/k\in J=[\min f',\max f']$ there is clearly exactly one stationary point $x_0$ such that $f'(x_0)=j/k$.

To this end, with the purpose of locating the stationary points for all $j/k\in J$, we define the ``dual curve" $f^*$ of $f$ as
\begin{equation}\label{e2.12}
f^*(y)=yh(y)-f(h(y)),
\end{equation}
where $h$ is the inverse function of $f'$. Note that 
\begin{equation}\label{e2.13}
(f^*)'(y)=h(y)+yh'(y)-f'(h(y))h'(y)=h(y)
\end{equation}
and hence 
\begin{equation}\label{e2.10}
(f^*)''(y)=h'(y)=\frac1{f''(h(y))}
\end{equation}
and
\begin{equation}\label{e2.11}
\frac1{c_2}\le |(f^*)''(y)|\le \frac1{c_1}.
\end{equation}
The operation $*$ , as it stands, is indeed an involution, meaning that when applying it twice in a row one gets back to the function which one starts with. The same construction, in one form or another, seems to appear in many areas of mathematics, like number theory, harmonic analysis and thermodynamics. It is also known as the Legendre transform in thermodynamics.  

Notice that $f^*$ is defined on $J$. For fixed $j$ and $k$, the stationary point is $x_0=h(j/k)$. Hence by the mean value theorem and \eqref{e2.11}
$$\min\left\{(f^*)'\left(\frac{j}k\right)-\eta,\xi-(f^*)'\left(\frac{j}k\right)\right\}\ge \frac1{c_2}D\left(\frac{j}k\right)$$
where
$$D(y)=\min\left\{y-\min f',\max f'-y\right\}.$$
Now by \eqref{e2.12}, \eqref{e2.13} and \eqref{e2.10}, we have
$$kf(x_0)-jx_0=k\left(f\left(h\left(\frac{j}k\right)\right)-\frac{j}k h\left(\frac{j}k\right)\right)=-kf^*\left(\frac{j}k\right)$$
and
$$f''(x_0)=\frac1{(f^*)''\left(\frac{j}k\right)}.$$
Thus Lemma \ref{l3} implies
\begin{align*}
&\int_{I}e(q(kf(x)-jx))dx\\
=&\sqrt{\frac{|(f^*)''(j/k)|}{qk}}e\left(-qkf^*\left(\frac{j}k\right)\pm\frac18\right)+O\left(\frac{{D(j/k)}^{-1}+\log (qk)}{qk}\right).
\end{align*}
Again, the sign of $1/8$ is chosen the same as that of $f''$, which does not change sign through out $I$. Hence this yields
\begin{align}
&\sum_{q\le Q}q\int_{I}e(q(kf(x)-jx))dx\nonumber\\
=&\sqrt{\frac{|(f^*)''(j/k))|}{k}}\sum_{q\le Q}{q}^{1/2}e\left(-qkf^*\left(\frac{j}k\right)\pm\frac18\right)+O\left(Q\frac{{D(j/k)}^{-1}+\log (Qk)}{k}\right).\label{e2.5}
\end{align}

We treat the error first. 
$$\sum_{j\in\tilde{{J}_k}}Q\frac{{D(j/k)}^{-1}+\log (Qk)}{k}\ll Q(\log k+\log Q)\ll Q\log Q.$$
Thus its contribution to \eqref{e2.0} is 
\begin{equation}\label{e2.15}
\ll\sum_{k\le K}\left(\frac1{K+1}+\min\left(\delta,\frac1{ k}\right)\right)Q\log Q\ll Q(\log Q)^2.
\end{equation}
Now we observe that when $\|kf^*(j/k)\|\le Q^{-1}$, the phase $e(-qkf^*(j/k))$ does not oscillate much when $q\le Q$. So we treat this case trivially. On the other hand, when $\|kf^*(j/k)\|> Q^{-1}$, there are indeed a lot of cancellations that we can exploit on summing over $q$. So we naturally divide the sum on the right side of \eqref{e2.5} into two cases and obtain by partial summation that

\begin{equation}
\sum_{q\le Q}q^{1/2}e\left(-qkf^*\left(\frac{j}k\right)\right)\ll\left\{
\begin{array}{lll}
Q^{3/2}& \text{if}& \|kf^*(j/k)\|\le Q^{-1}\label{e2.6}\\
Q^{1/2}\|kf^*(j/k)\|^{-1}&\text{if}&\|kf^*(j/k)\|>Q^{-1}.
\end{array}
\right.
\end{equation}

Here we state some consequences of Huxley's bound \eqref{e1.0}.

\begin{lem}\label{l4}
Let $F\in\mathcal{F}(J;c_2^{-1},c_1^{-1})$. 
Then for any $Q, K>1$
\begin{equation}\label{e2.7} 
\sum_{k\le K}\sum_{\substack{j/k\in J\\\|kF(j/k)\|>Q^{-1}}}k^{-1/2}\|kF(j/k)\|^{-1/2}\ll K^{3/2}+K^{1/2}(\log K) Q^{1/2},
\end{equation}
\begin{equation}\label{e2.8}
\sum_{k\le K}\sum_{\substack{j/k\in J\\\|kF(j/k)\|>Q^{-1}}}k^{-1/2}\|kF(j/k)\|^{-1}\ll_\varepsilon K^{3/2}Q^\varepsilon+K^{1/2}(\log K) Q,
\end{equation}
\begin{equation}\label{e2.9}
\sum_{k\le K}\sum_{\substack{j/k\in J\\\|kF(j/k)\|\le Q^{-1}}}k^{-1/2}\ll_\varepsilon K^{3/2}Q^{\varepsilon-1}+K^{1/2}\log K,
\end{equation}

\begin{equation}\label{e2.17} 
\sum_{K_1<k\le K_2}\sum_{\substack{j/k\in J\\\|kF(j/k)\|>Q^{-1}}}k^{-3/2}\|kF(j/k)\|^{-1/2}\ll K_2^{1/2}+K_1^{-1/2}(\log K_1) Q^{1/2},
\end{equation}
\begin{equation}\label{e2.18}
\sum_{K_1<k\le K_2}\sum_{\substack{j/k\in J\\\|kF(j/k)\|>Q^{-1}}}k^{-3/2}\|kF(j/k)\|^{-1}\ll_\varepsilon K_2^{3/2}Q^\varepsilon+K_1^{-1/2}(\log K_1) Q,
\end{equation}
\begin{equation}\label{e2.19}
\sum_{K_1<k\le K_2}\sum_{\substack{j/k\in J\\\|kF(j/k)\|\le Q^{-1}}}k^{-3/2}\ll_\varepsilon K_2^{1/2}Q^{\varepsilon-1}+K_1^{-1/2}\log K_1.
\end{equation}

\end{lem}

\eqref{e2.7}, \eqref{e2.8} and \eqref{e2.9} are essentially the same as Lemma 2.3 in \cite{VV}, while \eqref{e2.17}, \eqref{e2.18} and \eqref{e2.19} can be also easily derived by partial summation from \eqref{e1.0}.

Note that $j\in \tilde{J}_k$ implies that $j/k\in J$. Now taking $F=f^*$ in Lemma \ref{l4} and by partial summation we obtain from \eqref{e2.6} that
\begin{align*}
&\sum_{k\le K}b_k\left|\sum_{j\in\tilde{J}_k}\sqrt{\frac{|(f^*)''(j/k)|}{k}}\sum_{q\le Q}q^{1/2}e\left(-qkf^*\left(\frac{j}k\right)\pm\frac18\right)\right|\\\ll & Q^{3/2}\sum_{k\le K}b_k\sum_{\substack{j/k\in J\\\|kf^*(j/k)\|\le Q^{-1}}}k^{-1/2}+ Q^{1/2}\sum_{k\le K}b_k\sum_{\substack{j/k\in J\\\|kf^*(j/k)\|>Q^{-1}}}k^{-1/2}\|kf^*(j/k)\|^{-1}\\
\ll &\Sigma_1+\Sigma_2+\Sigma_3 
\end{align*}
where
$$\Sigma_1=\frac{Q^{3/2}}K\sum_{k\le K}\sum_{\substack{j/k\in J\\\|kf^*(j/k)\|\le Q^{-1}}}k^{-1/2}+ \frac{Q^{1/2}}K\sum_{k\le K}\sum_{\substack{j/k\in J\\\|kf^*(j/k)\|>Q^{-1}}}k^{-1/2}\|kf^*(j/k)\|^{-1},$$
$$\Sigma_2=\delta{Q^{3/2}}\sum_{k\le \frac1\delta}\sum_{\substack{j/k\in J\\\|kf^*(j/k)\|\le Q^{-1}}}k^{-1/2}+ \delta{Q^{1/2}}\sum_{k\le \frac1\delta}\sum_{\substack{j/k\in J\\\|kf^*(j/k)\|>Q^{-1}}}k^{-1/2}\|kf^*(j/k)\|^{-1}$$
and
$$\Sigma_3={Q^{3/2}}\sum_{\frac1\delta<k\le K}\sum_{\substack{j/k\in J\\\|kf^*(j/k)\|\le Q^{-1}}}k^{-3/2}+ {Q^{1/2}}\sum_{\frac1\delta<k\le K}\sum_{\substack{j/k\in J\\\|kf^*(j/k)\|>Q^{-1}}}k^{-3/2}\|kf^*(j/k)\|^{-1}.$$
Now by \eqref{e2.8} and \eqref{e2.9} we have
$$\Sigma_1\ll_\varepsilon K^{1/2}Q^{1/2+\varepsilon}+K^{-1/2}Q^{3/2}\log K$$
and
$$\Sigma_2\ll_\varepsilon \delta^{-1/2}Q^{1/2+\varepsilon}+\delta^{1/2} Q^{3/2}\log\frac1\delta.$$
And by \eqref{e2.18} and \eqref{e2.19} we have
$$\Sigma\ll_\varepsilon  K^{1/2}Q^{1/2+\varepsilon}+\delta^{1/2} Q^{3/2}\log \frac1\delta.$$
Recall that
$$K=Q\quad\text{and}\quad \delta\ge Q^{-1}$$
so
$$\Sigma_1+\Sigma_2+\Sigma_3 \ll_\varepsilon \delta^{\frac12} Q^{\frac32}\log\frac1\delta+Q^{1+\varepsilon}.$$
Hence combining the above estimates with \eqref{e2.0}, \eqref{e2.14} and \eqref{e2.15} yields
$$N_f(Q,\delta)-|I|\delta Q^2\ll_\varepsilon \delta^{\frac12} Q^{\frac32}\log\frac1\delta+Q^{1+\varepsilon} $$
for $\delta\ge Q^{-1}$. Finally note that $N_f(Q,\delta)$ is an increasing function in $\delta$ for fixed $Q$. Hence when $\delta< Q^{-1}$
$$N_f(Q,\delta)\le N_f(Q,Q^{-1})\ll_\varepsilon Q^{1+\varepsilon}.$$
In any case, Theorem \ref{t1} holds.

\section{Weakly non-degenerate planar curves}\label{s3}
In \S\ref{s2} we apply the method of stationary phase (Lemma \ref{l3}) to the oscillatory integral \eqref{e2.3}. Lemma \ref{l3} requires that the phase function be $C^2$ and its second derivative be Lipschitz continuous. In this section, we prove Theorem \ref{t2}, a weaker version of Theorem \ref{t1} under no Lipschitz assumption. Theorem \ref{t2} is nevertheless general enough to extend Khintchine type theorems to the class of weakly non-degenerate planar curves, which is first introduced in \cite{BZ}.   

Our strategy is to show Theorem \ref{t2} holds for $\mathcal{F}(I;c_1,c_2)$ first. It is readily verified that the $O$-constants in our argument  are uniform for all $f\in\mathcal{F}(I;c_1,c_2)$, hence a simple limiting argument extends the same result to the topological closure $\overline{\mathcal{F}}(I;c_1,c_2)$ of $\mathcal{F}(I;c_1,c_2)$. Thus without loss of generality, we can assume from now on $f\in\mathcal{F}(I;c_1,c_2)$.

For the time being, we assume 
$$\delta\ge Q^{-\frac23}$$
and
$$K=Q^{\frac23}.$$

We follow the argument in \S\ref{s2} up to \eqref{e2.4}, right before we apply Lemma \ref{l3}. We summarize it below.

\begin{equation}
N_f(Q,\delta)-|I|\delta Q^2\nonumber
\ll\sum_{k=1}^{K}b_k\left|\sum_{q\le Q}q\sum_{j\in\tilde{J}_k}\int_{I}e(q(kf(x)-jx))dx\right|+E\label{e3.1}
\end{equation}
where
$$b_k=\frac1{K+1}+\min\left(2\delta,\frac1{\pi k}\right)$$
and 
\begin{align}
E&\ll\frac{Q^2}{K}+\delta Q+Q(\log K)^2+K^{-\frac12}Q^{\frac32}+\delta^{\frac12}Q^{\frac32}\nonumber\\
&\ll Q^{\frac43}+\delta^{\frac12}Q^{\frac32}\label{e3.1}.
\end{align}

Although we cannot officially apply Lemma \ref{l3}, we still bound the contribution to \eqref{e2.3} from those points which are near or away from the critical point separately. In particular, we can still define the dual curve $f^*(y)$ as in \eqref{e2.12}. The relations between $f$ and $f^*$ still hold.

By \eqref{e2.20} we conclude the terms with $\|kf^*(j/k)\|\le Q^{-1}$ contribute
\begin{align}
\ll&\sum_{k=1}^{K}b_k\sum_{q\le Q}q\sum_{\substack{j/k\in{J}\\\|kf^*(j/k)\|\le Q^{-1}}}\frac1{\sqrt{qk}}\nonumber\\
\overset{\eqref{e2.9}\&\eqref{e2.19}}{\ll_\varepsilon}&K^{1/2}Q^{1/2+\varepsilon}+K^{-1/2}(\log K) Q^{3/2}+\delta^{-1/2}Q^{1/2+\varepsilon}+\delta^{1/2}\left(\log\frac1\delta\right) Q^{3/2}\nonumber\\
\ll&\delta^{1/2}\left(\log\frac1\delta\right) Q^{3/2}+Q^{4/3}\label{e3.2}.
\end{align}

Now we turn to the terms with $\|kf^*\left(\frac{j}k\right)\|>Q^{-1}$.

For fixed $j$ and $k$ with $j/k\in J$, recall that $x_0$ denotes the critical point with $f'(x_0)=j/k$. For simplicity, we denote
$$\beta=\|kf(x_0)-jx_0\|=\|kf^*\left(\frac{j}k\right)\|.$$
 let
$$\mathcal{A}(j,k)=\left\{x\in I:|x-x_0|>\sqrt{\frac{\beta}{c_2k}}\right\}.$$
The Lagrange mean value theorem implies that for any $x\in\mathcal{A}(j,k)$ there exists $x_1\in I$ such that
$$kf'(x)-j=(x-x_0)kf''(x_1).$$
Thus
$$|kf'(x)-j|\gg \sqrt{k\beta}.$$
Hence by integration by parts (or instead Lemma 1 in \cite[Chapter 3]{Mo}), we obtain
\begin{equation}\label{e3.3}
\sum_{q\le Q}\int_{\mathcal{A}(j,k)}e(q(kf(x)-jx))dx\ll \frac{Q}{\sqrt{k\beta}}.
\end{equation}

On the other hand, let
$$\mathcal{B}(j,k)=I\backslash\mathcal{A}(j,k).$$
By Taylor's theorem, for any $x\in\mathcal{B}(j,k)$, there exists $x_2\in I$ such that
$$
kf(x)-jx-(kf(x_0)-jx_0)=\frac12kf''(x_2)(x-x_0)^2
$$
the right side of which is bounded by
$$\frac12c_2k\frac{\beta}{c_2k}=\frac12\beta.$$
Consequently for any $x\in\mathcal{B}(j,k)$,
$$\|kf(x)-jx\|\asymp \beta.$$
Therefore
\begin{align}
\int_{\mathcal{B}(j,k)}\sum_{q\le Q}qe(q(kf(x)-jx))dx
\ll& Q\beta^{-1}|\mathcal{B}(j,k)|\nonumber\\
\ll& \frac{Q}{\sqrt{k\beta}}.\label{e3.4}
\end{align}
Merging \eqref{e3.3} and \eqref{e3.4} reveals
$$\sum_{q\le Q}q\int_Ie(q(kf(x)-jx))dx\ll Q\|kf^*(j/k)\|^{-1/2}k^{-1/2}.$$
Hitherto, by \eqref{e2.7} and \eqref{e2.17}
\begin{align*}
&\sum_{k=1}^{K}b_k\sum_{q\le Q}q\sum_{\substack{j/k\in{J}\\\|kf^*(j/k)\|> Q^{-1}}}\|kf^*(j/k)\|^{-1/2}k^{-1/2}\\
\ll&K^{1/2}Q+K^{-1/2}(\log K)Q^{3/2}+\delta^{-1/2}Q+\delta^{1/2}\left(\log\frac1\delta\right)Q^{3/2}\\
\ll&\delta^{1/2}\left(\log\frac1\delta\right)Q^{3/2}+Q^{4/3}
\end{align*}
on recalling our assumption $K=Q^{2/3}$ and $\delta\ge Q^{-2/3}$.

Now combining this with various errors \eqref{e3.1} and \eqref{e3.2} we have obtained so far, we conclude that for any $f\in\mathcal{F}(I;c_1,c_2)$ and any $\delta\ge Q^{-2/3}$
$$
N_f(Q,\delta)-|I|\delta Q^2\ll \delta^{1/2}\left(\log\frac1\delta\right)Q^{3/2}+Q^{4/3}.
$$
Again notice that for fixed $f$ and $Q$, by definition $N_f(Q,\delta)$ is an increasing function in $\delta$. So when $\delta<Q^{-2/3}$
$$N_f(Q,\delta)\le N_f(Q,Q^{-2/3})\ll Q^{7/6}\log Q+Q^{4/3}\ll Q^{4/3}.$$
So, we have proven Theorem \ref{t2}.

\section{The method of stationary phase}\label{s4}
In this section, we prove Lemma \ref{l3}. The proof given here is essentially embedded in \cite[\S 4]{VV}. We believe that Lemma \ref{l3} as stated in this paper is of independent interest and may have further applications. In standard reference books like \cite{Mo,Pi,St}, it is always assumed that the phase function is at least $C^4$ or even $C^\infty$. The advantage of Lemma \ref{l3} is that it only requires $f$ be $C^2$ and $f''$ be Lipschitz continuous.\footnote{We note that in \cite[Theorem 2.7.1]{Pi} a conclusion of a comparable quality to that of Lemma \ref{l3} is claimed without the latter condition. Unfortunately we observe a serious oversight in the proof, which requires additional differentiability condition on the phase function to rectify.}

According to the definition of $\mathcal{F}(I;c_1,c_2)$, we know $\phi''(x)$ does not change sign on $I$. So without loss of generality we may assume $\phi''(x)>0$. The other case is formally the same as taking the complex conjugate.

We extend the definition of $\phi$ to $\mathbb{R}$ by letting 
$$\phi(x)=\phi(b)+\phi'(b)(x-b)+\frac{\phi''(b)}{2}(x-b)^2, \text{   when } x>b$$
and
$$\phi(x)=\phi(a)+\phi'(a)(x-a)+\frac{\phi''(a)}{2}(x-a)^2, \text{   when } x<a.$$
Now notice that $\phi''$ is Lipschitz on $\mathbb{R}$ and still satisfies $c_1\le \phi''(x)\le c_2$. We shift the integral to an interval centered at $x=x_0$. Let 
$$\mu=\frac{b-a}2.$$
We consider the error
\begin{equation*}
\left(\int_a^b-\int_{x_0-\mu}^{x_0+\mu}\right)e(\lambda\phi(x))dx.
\end{equation*}
If $x\not\in[a,b]\cap[x_0-\mu,x_0+\mu]$, then $|x-x_0|>\kappa$, where $\kappa=\min\{x_0-a,b-x_0\}$. In this case, by the mean value theorem
$$|\phi'(x)|\ge c_1|x-x_0|>c_1\kappa.$$ Now applying for instance Lemma 1 in \cite[Chap. 3]{Mo} we know 
\begin{equation}\label{e4.1}
\left(\int_a^b-\int_{x_0-\mu}^{x_0+\mu}\right)e(\lambda\phi(x))dx\ll
 \frac1{\kappa\lambda}.
\end{equation}
Next we write
$$v(x)=\phi(x+x_0)-\phi(x_0).$$ Then
$$v(0)=v'(0)=0, \quad v''(0)=\phi''(x_0)$$
and
\begin{equation}\label{e4.0}
\int_{x_0-\mu}^{x_0+\mu}e(\lambda\phi(x))dx=e(\lambda\phi(x_0))\int_{-\mu}^{\mu}e(\lambda v(x))dx.
\end{equation}

Since $v''(x)=\phi''(x+x_0)>0$ and $v'(0)=0$, $v'$ is strictly positive on $(0,\infty)$. Thus $v$ is strictly increasing on $[0,\infty]$. Let $x(v)$ be the inverse function of $v(x)$ on $[0,\infty)$. Clearly $x(v)$ is twice continuously differentiable on $(0,\infty)$. For any $\nu\in(0,\mu)$ we have
\begin{equation}\label{e4.2}
\int_{\nu}^{\mu}e(\lambda v(x))dx=\int_{v(\nu)}^{v(\mu)}e(\lambda v)x'(v)dv.
\end{equation}
By the Lipschitz condition and Taylor's formula, we have for $x>0$
$$v''(x)=v''(0)+O(x)$$
and
$$v'(x)=v''(0)x+O(x^2)$$
and
$$v(x)=\frac12v''(0)x^2+O(x^3).$$
Furthermore for $0<x\le\mu$ we have
$$x(v)=\sqrt{\frac{2v}{v''(0)}}(1+O(v^{1/2})).$$
Therefore,
\begin{equation}\label{e4.3}
x'(v)=\frac1{v'(x(v))}=(2v''(0)v)^{-1/2}+O(1)
\end{equation}
and
\begin{equation}\label{e4.4}
x''(v)=-v''(x(v))(x'(v))^3=(8v''(0)v^3)^{-1/2}+O(v^{-1}).
\end{equation}

Now by integration by parts, the integral on the right hand of \eqref{e4.2} is
$$\left.\frac{e(\lambda v)x'(v)}{2\pi i\lambda}\right|_{v(\nu)}^{v(\mu)}-\int_{v(\nu)}^{v(\mu)}\frac{e(\lambda v)}{2\pi i \lambda}x''(v)dv.$$
The first term, by the approximation \eqref{e4.3}  is
\begin{align*}
&-\frac{e(\lambda v(\nu))x'(v(\nu))}{2\pi i\lambda}+O(\lambda^{-1})\\
=&-\frac{e(\lambda v(\nu))(2v''(0)v(\nu))^{-1/2}}{2\pi i\lambda}+O(\lambda^{-1})
\end{align*}
and the second term, by the approximation \eqref{e4.4} is
\begin{align*}
\int_{v(\nu)}^\infty\frac{e(\lambda v)(8v''(0)v^3)^{-1/2}}{2\pi i \lambda}dv+O\left(\frac{|\log\nu|+1}{\lambda}\right).
\end{align*} 

Now on combining the above estimates we have
\begin{align}
&\phantom{xxxxxxxxxxxx}\int_{\nu}^{\mu}e(\lambda v(x))dx\nonumber\\
=&-\frac{e(\lambda v(\nu))(2v''(0)v(\nu))^{-1/2}}{2\pi i\lambda}-\int_{v(\nu)}^\infty\frac{e(\lambda v)(8v''(0)v^3)^{-1/2}}{2\pi i \lambda}dv+O\left(\frac{|\log\nu|+1}{\lambda}\right)\nonumber\\
=&\int_{v(\nu)}^\infty\frac{e(\lambda v)}{\sqrt{2v''(0)v}}dv+O\left(\frac{|\log\nu|+1}{\lambda}\right).\label{e4.5}
\end{align}

On the other hand,
\begin{align}
\int_0^{\nu}e(\lambda v(x))dx&=\int_0^{\nu}e\left(\lambda \frac12v''(0)x^2\right)+O\left(\int_0^\nu\lambda x^3dx\right)\nonumber\\
&=\int_0^{\frac12v''(0)\nu^2}\frac{e(\lambda v)}{\sqrt{2v''(0)v}}dv+O(\lambda\nu^4)\nonumber\\
&=\int_0^{v(\nu)}\frac{e(\lambda v)}{\sqrt{2v''(0)v}}dv+O(\lambda\nu^4+\nu^2).\label{e4.6}
\end{align}
The last step in the above equation is justified because
$$\int_{\frac12v''(0)\nu^2}^{v(\nu)}\frac{e(\lambda v)}{\sqrt{2v''(0)v}}dv\ll \frac1{\sqrt{\nu^2}}\nu^3\ll \nu^2.$$

Hitherto, we obtain by adding up \eqref{e4.5} and \eqref{e4.6} that
$$\int_0^{\mu}e(\lambda v(x))dx=\int_0^{\infty}\frac{e(\lambda v)}{\sqrt{2v''(0)v}}dv+O\left(\frac{|\log\nu|+1}{\lambda}+\lambda\nu^4+\nu^2\right)$$
where the error is
$$O\left(\frac{\log\lambda}{\lambda}\right)$$
on choosing
$$\nu=\frac{\mu}{\sqrt{\lambda}}.$$
Hence performing the change of variable $v=\frac{z^2}{\lambda}$ gives
\begin{equation*}
\int_0^{\infty}\frac{e(\lambda v)}{\sqrt{2v''(0)v}}dv=\frac2{\sqrt{2v''(0)\lambda}}\int_0^\infty e(z^2)dz,
\end{equation*}
where the Fresnel integral 
$$\int_0^\infty e(z^2)dz=\frac{e^{i\frac{\pi}4}}{2\sqrt{2}}$$
by the standard contour integration method. We thus have
\begin{equation}\label{e4.7}
\int_0^{\mu}e(\lambda v(x))dx=\frac{e\left(\frac18\right)}{2\sqrt{v''(0)\lambda}}+O\left(\frac{\log\lambda}{\lambda}\right).
\end{equation}
The integral $\int_{-\mu}^{0}e(\lambda v(x))dx$ can be treated in exactly the same manner and shown to have the same contribution as \eqref{e4.7}. Now Lemma \ref{l3} follows from \eqref{e4.1}, \eqref{e4.0} and \eqref{e4.7}.

\section{A problem of Mazur}\label{s5}
Before the recent work in \cite{BDV}, little had been known about the existence of rational points near a planar curve for small $\delta$ say $\delta<Q^{-1/2}$. Motivated by Elkies \cite{El}, an explicit question of this type was raised by Barry Mazur \cite{Ma} who asks: ``\emph{given a smooth curve in the plane, how near to it can a point with rational coordinates get and still miss?}" Of course one wants to measure the distance to the curve of a rational point in terms of its denominator. So Mazur's question is concerned with those rational points counted in $\hat{N}_f(Q,\delta)$ that do not lie on the curve. There is a very rich theory about rational points on algebraic curves. Curves with similar geometric properties may exhibit quite different arithmetic behaviors. For instance, the circle $x^2+y^2=1$ centered at the origin with radius 1 certainly contains infinitely many rational points. But the larger circle $x^2+y^2=3$ has no rational points. A more delicate example is the Fermat curve $x^3+y^3=1$, which can be shown to have only two rational points $(1,0)$, $(0,1)$.

Now look at our Theorem \ref{t1}. It has a dominating main term $|I|\delta Q^2$ when $\delta>Q^{-1+o(1)}$, which has nothing to do with the arithmetic of the curve. So one may wonder where does the arithmetic go? It is in the error term! Let $\delta$ tend to 0. One sees that the term $Q^{1+\varepsilon}$ does not change even though the main term goes to 0. Actually there is a good explanation for this. Rational points on the curve are always counted in $\hat{N}_f(Q,\delta)$ regardless of the value of $\delta$. So the main term corresponds to those rational points which are genuinely near the curve, while the term $Q^{1+\varepsilon}$ corresponds to rational points lying on the curve. We do not seem to have a good explanation for the middle term $\delta^{1/2}(\log \delta^{-1})Q^{3/2}$ though we do observe that it is essentially the geometric mean of $\delta Q^2$ and $Q^{1+\varepsilon}$ hence is always bounded by one of them. So Theorem \ref{t1} produces a general answer to Mazur's question; namely, for any $C^3$ curve with curvature bounded from above and below by positive numbers, there always exists a rational point with denominator at most $Q$ that is within distance   $O(Q^{-2+o(1)})$ from the curve. Moreover, the quantity of these rational points satisfies an asymptotic formula that matches the probabilistic heuristics (see the discussion below Theorem \ref{t1}).

On the other hand, if $\delta=o(Q^{-1})$, then we have a completely different story. One in general does not expect to find any rational points near the curve except for those already on the curve. Indeed, it is shown in \cite{BD} that when $\delta=o(Q^{-1})$ the rational points in question cannot miss the curve if it is rational quadratic. This  shows that the answer above is in general best possible. It is also shown in \cite{BD} that when $\delta=o(Q^{-k+1})$ the rational points in question cannot miss the Fermat curve $x^k+y^k=1$ when $k\ge3$. But we now know there are at most 4 rational points on it thanks to Wiles \cite{Wi}! So we conclude that $\hat{N}_f(Q,\delta)\le 4$ when $f$ is the function representing the Fermat curve $x^k+y^k=1$ and $\delta=o(Q^{-k+1})$. The upshot here is that when $\delta$ is small, the arithmetic of the curve starts to kick in; therefore the counting problem considered in this article is equivalent to counting points on the curve, which is commonly known to be very delicate.

It is a very interesting phenomenon that when $\delta$ is not too small, $\hat{N}_f(Q,\delta)$ is governed by the ``coarse" geometric property say the curvature and when $\delta$ is small, the problem is governed by the ``fine" arithmetic property. To be more precise, in the former case rational points near the curve do exist and are evenly distributed (they behave exactly as expected), while in the latter case, we believe they must lie on the curve and so their quantity must depend on the arithmetic nature of the curve and thus fluctuates a lot from curve to curve. This discussion, if making sense at all, should suggest that there exists a turning point for each curve (at least for algebraic curves). We see above that, at least for rational quadratic curves, the turning point occurs when $\delta\approx Q^{-1}$. We also know that the turning point for the Fermat curve $x^k+y^k=1$ with $k\ge 3$ occurs for some $\delta$ between $Q^{-k+1}$ and $Q^{-1}$, though determining this exactly is well beyond the reach of any current machinery. In general, it should be very difficult to determine the exact threshold for any generic curve. One sees that this is actually in the same spirit as Mazur's \emph{near-misses problem} \cite{Ma}.  There is a similar discussion about Mazur's problem in \cite{Wa}.

\section{Uniform distribution}\label{s6}
Looking at \eqref{e2.1}, we may recall the theory of uniform distribution. Indeed, what we really have established in Theorem \ref{t1} can be naturally regarded as a type of strong uniform distribution result. We recall some basic definitions first.

Let $\{u_n\}$, $n=1,2,3,\cdots$ be a sequence of real numbers. Let $Z(N;\alpha,\beta)$ count the number of $n\le N$ for which $u_n\in(\alpha,\beta)\pmod{1}$ with $\alpha<\beta<\alpha+1$.  The discrepancy function of the sequence $\{u_n\}$ is defined by

$$D(N):=\sup_{\alpha,\beta}\left|D(N;\alpha,\beta)\right|$$
where
$$
D(N;\alpha,\beta):=Z(N;\alpha,\beta)-(\beta-\alpha)N.
$$

The sequence $\{u_n\}$ is said to be \emph{uniformly distributed} if
$$\lim_{N\rightarrow\infty}\frac{D(N)}N=0.$$

If we take $\{u_n\}$ to be the sequence $\{qf(a/q)\}$, with $a/q\in I$ and $q\le Q$ as $Q\rightarrow\infty$, for some real valued function $f:I\rightarrow\mathbb{R}$. Let $N=N(Q)$ be the number of terms in this sequence with $q\le Q$. Thus $N\asymp |I|Q^2/2$. We use $\{x\}$ to denote the fractional part of $x$. Then by making a very minor change to the proof of Theorem \ref{t1}, we can prove the following theorem.

\begin{thm}[Strong uniform distribution]\label{t5} 
Let $f$ satisfy the same conditions with Theorem \ref{t1} and  the sequence $\{u_n\}$ be given above. Then for any $\varepsilon>0$,

$$D(N;\alpha,\beta)\ll (\beta-\alpha)^\frac12\log(\beta-\alpha)^{-1}N^\frac34+N^{\frac12+\varepsilon}$$
and
$$D(N)\ll N^\frac34.$$

\end{thm}

\section{Further discussions on general manifolds}\label{s7}

It remains interesting to see whether the method employed in this paper can be extended to establish counting results of the same type for submanifolds in Euclidean space of higher dimension. A lower bound of the expected order of magnitude was established by Beresnevich \cite{Be} but the upper bound in general remains a major open question in this area. The framework we have laid here should extend to general manifolds. But we crucially need Huxley's bound \eqref{e1.0}. In some sense, our Theorem \ref{t1} is equivalent to Huxley's bound, even though the former looks stronger than the latter. The secret here is that the harmonic analysis lets one transfer back and forth between the physical space and the frequency space and what it does is just to send a bound in one space to another equivalent form in the dual space. Huxley's bound is used exactly as an input in our argument. And indeed, if checking our argument carefully, one observes that the $\varepsilon$ in  $\delta^{1-\varepsilon}Q^2$ of Huxley's bound is exactly the same $\varepsilon$ that appears in $Q^{1+\varepsilon}$ in the error of Theorem \ref{t1}. So the use of Huxley is the real bottleneck of our argument. If one wants to prove anything as sharp as Theorem \ref{t1} about general manifolds in higher dimensions along the same line, one has to either generalize Huxley to higher dimensions or produce an independent method that does the same, though neither seems easy. However, it seems likely that weaker upper bounds with shorter range of $\delta$ can be established for some manifolds with non-vanishing principal curvatures. We will return to this in a future publication.

\subsection*{Acknowledgments}

The author would like to thank Prof. John Friedlander for constant help and encouragement during the preparation of this article.

\end{document}